\documentclass{amsart}
\usepackage{epsfig}
\usepackage{amssymb}
\usepackage{latexsym}

\def\ll{\lambda}
\def\strict{\mathrm{strict}}
\def\weak{\mathrm{weak}}
\def\Qsym{\mathcal Qsym}
\def\Q{\mathbb Q}
\def\C{\mathcal C}
\def\bd{\mathrm{bd}}
\def\P{\mathbb P}
\def\N{\mathbb N}
\def\T{\mathbb T}
\def\Z{\mathbb Z}
\def\wt{\mathrm{wt}}

\newtheorem{theorem}{Theorem}[section]
\newtheorem{lemma}[theorem]{Lemma}
\newtheorem{proposition}[theorem]{Proposition}
\newtheorem{example}[theorem]{Example}

\newtheorem{conjecture}[theorem]{Conjecture}
\theoremstyle{definition}
\newtheorem{definition}[theorem]{Definition}

\title{Cell transfer and monomial positivity}
\author{Thomas Lam and Pavlo Pylyavskyy}
\address{Department of Mathematics, Harvard University, Cambridge, MA, 02138}
\email{tfylam@math.harvard.edu}
\thanks{T.L. was supported in part by NSF
DMS-0600677.}
\address{Department of Mathematics, M.I.T., Cambridge, MA, 02139}
\email{pasha@mit.edu}

\begin{document}

\begin{abstract}
We give combinatorial proofs that certain families of differences
of products of Schur functions are monomial-positive.  We show in
addition that such monomial-positivity is to be expected of a
large class of generating functions with combinatorial definitions
similar to Schur functions.  These generating functions are
defined on posets with labelled Hasse diagrams and include for
example generating functions of Stanley's $(P,\omega)$-partitions.
\end{abstract}
\maketitle

\section{Introduction}
\label{sec:introduction}



The Schur functions $s_\ll$ form a basis of the ring of symmetric functions $\Lambda$.
They have a remarkable number of combinatorial and algebraic
properties, and are simultaneously the irreducible characters of
$GL(N)$ and representatives of Schubert classes in the cohomology
$H^*(Gr_{kn})$ of the Grassmannian; see~\cite{Mac,EC2}.  In recent
years, a lot of work has gone into studying whether certain
expressions of the form
\begin{equation}
\label{eq:prod} s_\ll s_\mu - s_\nu s_\rho
\end{equation}
are expressible as a non-negative linear combination of Schur
functions.  See for example~\cite{BM,BBR,FFLP,Oko}.

The first aim of this article is to provide a large class of
expressions of the form (\ref{eq:prod}) which are {\it
monomial-positive}, that is, expressible as a non-negative linear
combination of monomial symmetric functions.  In particular, we show that
(\ref{eq:prod}) is monomial-positive when $\ll =\nu \vee \rho$ and
$\mu = \nu \wedge \rho$ are the union and intersections of the
Young diagrams of $\nu$ and $\rho$.  However, we show in addition
that such monomial-positivity is to be expected of many families
of generating functions with combinatorial definitions similar to
Schur functions, which are generating functions for semistandard
Young tableaux.

We define a new combinatorial object called a {\it $\T$-labelled
poset} and given a $\T$-labelled poset $(P,O)$ we define another
combinatorial object which we call a {\it $(P,O)$-tableau}.  These
$(P,O)$-tableaux include as special cases standard Young tableaux,
semistandard Young tableaux, shifted tablueax, cylindric tableaux,
plane partitions, and Stanley's $(P,\omega)$-partitions.  Our main
theorem is the {\it cell transfer theorem}.   It says that for a
fixed $\T$-labelled poset $(P,O)$, one obtains many expressions of
the form (\ref{eq:prod}) which are monomial-positive, where the
Schur functions in (\ref{eq:prod}) are replaced by generating
functions for $(P,O)$-tableaux.

We conjecture that our cell-transfer results for Schur functions
hold not just for monomial-positivity but also for Schur-positivity.
This conjecture is proved in~\cite{LPP}.  In another direction, we
strengthen the results of the present article in the case of
generating functions of $(P,\omega)$-partitions in~\cite{LP}. In
this case, cell transfer is positive in terms of fundamental
quasisymmetric functions.

\medskip
{\bf Acknowledgements.} We would like to thank our advisor Richard
Stanley, for interesting conversations concerning this problem.  We thank the
anonymous referee for many helpful suggestions.

\section{Posets and Tableaux}

Let $(P, \leq)$ be a possibly infinite poset.  Let $s, t \in P$.
We say that $s$ {\it covers} $t$ and write $s \gtrdot t$ if for
any $r \in P$ such that $s \geq r \geq t$ we have $r = s$ or $r =
t$. The {\it Hasse diagram} of a poset $P$ is the graph with
vertex set equal to the elements of $P$ and edge set equal to the
set of covering relations in $P$.  If $Q \subset P$ is a subset of
the elements of $P$ then $Q$ has a natural induced subposet
structure. If $s, t \in Q$ then $s \leq t$ in $Q$ if and only if
$s \leq t$ in $P$. Call a subset $Q \subset P$ {\it connected} if
the elements in $Q$ induce a connected subgraph in the Hasse
diagram of $P$.

An {\it order ideal} $I$ of $P$ is an induced subposet of $P$ such
that if $s \in I$ and $s \geq t \in P$ then $t \in I$.  A subposet
$Q \subset P$ is called {\it convex} if for any $s, t \in Q$ and
$r \in P$ satisfying $s \leq r \leq t$ we have $r \in Q$.
Alternatively, a convex subposet is one which is closed under
taking intervals.  A convex subset $Q$ is determined by specifying
two order ideals $J$ and $I$ so that $J \subset I$ and $Q = \{s
\in I \mid s \notin J\}$. We write $Q = I/J$.  If $s \notin Q$
then we write $s < Q$ if $s < t$ for some $t \in Q$ and similarly
for $s > Q$. If $s \in Q$ or $s$ is incomparable with all elements
in $Q$ we write $s \sim Q$. Thus for any $s \in P$, exactly one of
$s <Q$, $s > Q$ and $s \sim Q$ is true.

Let $\P$ denote the set of positive integers and $\Z$ denote the
set of integers.  Let $\T$ denote the set of all weakly increasing
functions $f: \P \rightarrow \Z \cup \{\infty\}$.

\begin{definition}
A {\it $\T$-labelling} $O$ of a poset $P$ is a map $O : \{(s,t)\in
P^2 \mid s \gtrdot t \} \rightarrow \T$ labelling each edge
$(s,t)$ of the Hasse diagram by a weakly increasing function
$O(s,t):\P \rightarrow \Z \cup \{\infty\}$.  A {\it $\T$-labelled
poset} is an an ordered pair $(P,O)$ where $P$ is a poset, and $O$
is a $\T$-labelling of $P$.
\end{definition}

We shall refer to a $\T$-labelled poset $(P,O)$ as $P$ when no
ambiguity arises.  If $Q \subset P$ is a convex subposet of $P$
then the covering relations of $Q$ are also covering relations in
$P$.  Thus a $\T$-labelling $O$ of $P$ naturally induces a
$\T$-labelling $O|_Q$ of $Q$.  We denote the resulting
$\T$-labelled poset by $(Q,O):=(Q,O|_Q)$.

\begin{definition}
A $(P, O)$-tableau is a map $\sigma: P \to \mathbb P$ such that
for each covering relation $s \lessdot t$ in $P$ we have
\[
\sigma(s) \leq O(s,t)(\sigma(t)).
\]
If $\sigma: P \to \mathbb P$ is any map, then we say that $\sigma$
{\it respects} $O$ if $\sigma$ is a $(P, O)$-tableau.
\end{definition}

\begin{figure}
\begin{center}
\epsfig{file=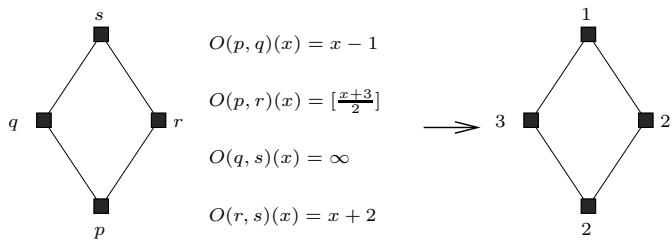}
\end{center}
\caption{An example of a $\T$-labelled
poset $(P,O)$ and a $(P,O)$-tableaux.}\label{fig:pos1}
\end{figure}

Figure~\ref{fig:pos1} contains an example of a $\T$-labelled poset
$(P, O)$ and a corresponding $(P,O)$-tableau.

Denote by $\mathcal{A}(P,O)$ the set of all $(P,O)$-tableaux. If
$P$ is finite then one can define the formal power series
$K_{P,O}(x_1,x_2,\ldots) \in \Q[[x_1,x_2,\ldots]]$ by
\[
K_{P,O}(x_1,x_2,\ldots) = \sum_{\sigma \in \mathcal{A}(P,O)}
x_1^{\# \sigma^{-1}(1)} x_2^{\# \sigma^{-1}(2)} \cdots.
\]
The composition $\wt(\sigma) = (\# \sigma^{-1}(1), \#
\sigma^{-1}(2), \ldots)$ is called the weight of $\sigma$.

Our $(P,O)$-tableaux can be viewed as a generalization of
Stanley's $(P,\omega)$-partitions and also of McNamara's oriented
posets; see~\cite{EC2,McN}.



\begin{example}
\label{ex:schur} Any Young diagram $P=\lambda$ can be considered as
a $\T$-labelled poset $(\lambda, O_\lambda)$.  The elements of
$\lambda$ are given by the squares of the Young diagram.  A square
$s \in \lambda$ is less than another square $s' \in \lambda$ if and
only if $s$ lies (weakly) above and to the left of $s'$ when
$\lambda$ is drawn in the English notation.  The covering relations
$(s \lessdot s')$ are given by pairs of squares sharing an edge.

The labelling $O_\lambda$ of the Hasse diagram of $\lambda$ is
obtained as follows.  An edge $(s \lessdot s')$ where $s$ and $s'$
lie in the same row is labelled with the function $f^{\weak}(x) = x$
and if $s$ and $s'$ lie in the same column, the edge is labelled
with the function $f^{\strict}(x) = x-1$.  A $(\lambda,
O_\lambda)$-tableau is then just a semistandard Young tableaux and
we have the equality $$K_{\lambda,O_\lambda}(x_1,x_2,\cdots)=
s_{\lambda}(x_1,x_2,\cdots),$$ where $s_{\lambda}(x_1,x_2,\cdots)$
is the Schur function labelled by $\lambda$ (see
Section~\ref{sec:sym}).

More generally, suppose $\lambda$ and $\mu$ are two partitions satisfying
$\mu \subset \lambda$.  The skew shape $\lambda/\mu$
can be considered a $\T$-labelled poset, and in this way we obtain
the skew Schur functions.
\end{example}

\begin{example}
\label{ex:cylindric} Another interesting example is given by
cylindric tableaux and cylindric Schur functions.  Let $1 \leq k <
n$ be two positive integers. Let $\C_{k,n}$ be the quotient of
$\Z^2$ given by
\[
\C_{k,n} = \Z^2/(k-n,k)/Z.
\]
In other words, the integer points $(a,b)$ and $(a+k-n,b+k)$ are
identified in $\C_{k,n}$.  We can give $\C_{k,n}$ the structure of
a poset by the generating relations $(i,j) \lessdot (i+1,j)$ and
$(i,j) \lessdot (i,j+1)$.  We give $\C_{k,n}$ a $\T$-labelling $O$
by labelling the edges $(i,j) \lessdot (i+1,j)$ with the function
$f^{\weak}(x) = x$ and the edges $(i,j) \lessdot (i,j+1)$ with the
function $f^{\strict}(x) = x-1$.  A finite convex subposet $P$ of
$\C_{k,n}$ is known as a {\it cylindric skew shape};
see~\cite{GK,Pos,McN}.  The $(P,O)$-tableau are known as {\it
semistandard cylindric tableaux} of shape $P$ and the generating
function $K_{P,O}(x_1,x_2,\cdots)$ is the cylindric Schur function
defined in~\cite{BS,Pos}.
\end{example}

\begin{example}
Let $N$ be the number of elements in a poset $P$, and let $\omega: P
\longrightarrow [N]$ be a bijective labelling of elements of $P$
with numbers from $1$ to $N$. Recall that a $(P,\omega)$-partition
(see~\cite{EC2}) is a map $\sigma: P \longrightarrow \mathbb P$ such
that $s \leq t$ in $P$ implies $\sigma(s) \leq \sigma(t)$, while if
in addition $\omega(s)>\omega(t)$ then $\sigma(s) < \sigma(t)$.
Label now each edge $(s,t)$ of the Hasse diagram of $P$ with
$f^{\weak}$ or $f^{\strict}$, depending on whether $\omega(s) \leq
\omega(t)$ or $\omega(s) > \omega(t)$ correspondingly. It is not
hard to see that for this labelling $O$ the $(P,O)$-tableaux are
exactly the $(P,\omega)$-partitions.  Similarly, if we allow any
labelling of the edges of $P$ with $f^{\weak}$ and $f^{\strict}$, we
get the {\it {oriented posets}} of McNamara; see~\cite{McN}.
\end{example}

\section{The Cell Transfer Theorem}
A generating function $f \in \Q[[x_1,x_2,\ldots]]$ is {\it
monomial-positive} if all coefficients in its expansion into
monomials are non-negative.  If $f$ is actually a symmetric
function then this is equivalent to $f$ being a non-negative
linear combination of monomial symmetric functions.

Let $(P,O)$ be a $\T$-labelled poset.  Let $Q$ and $R$ be two
finite convex subposets of $P$.  The subset $Q \cap R$ is also a
convex subposet.  Define two more subposets $Q \wedge R$ and $Q
\vee R$ by
\begin{equation}\label{eq:wedge}
Q \wedge R = \{s \in R \mid s < Q \} \cup \{s \in Q \mid s \sim R
\;\text{or}\; s < R\}
\end{equation}
and
\begin{equation}\label{eq:vee}
Q \vee R = \{s \in Q \mid s > R \} \cup \{s \in R \mid s \sim Q
\;\text{or}\; s > Q\}.
\end{equation}

Observe that the operations $\vee$, $\wedge$ are not commutative,
and that $Q \cap R$ is a convex subposet of both $Q \vee R$ and $Q
\wedge R$.

\begin{example}
\label{ex:chain}
Let $P_n$ denote the chain with $n+1$ elements labeled $\{0,1,\ldots,n\}$.  Then the convex subposets of
$P_n$ are the intervals $[i,j]$ where $0 \leq i \leq j \leq n$ and $[i,j]$ is isomorphic to the
chain with $j-i$ elements.  Let $Q = [i,j]$ and $R = [i',j']$ and assume that $i \leq i'$.  Then we have the following two cases:
\begin{enumerate}
\item If $j \leq j'$ then $Q\wedge R = Q$ and $Q \vee R = R$.
\item If $j \geq j'$ then $Q \wedge R = [i,j']$ and $Q \vee R = [i',j]$.
\end{enumerate}
This example leads to interesting combinatorics which we study
further in~\cite{LP}.
\end{example}

\begin{figure}
\begin{center}
\epsfig{file=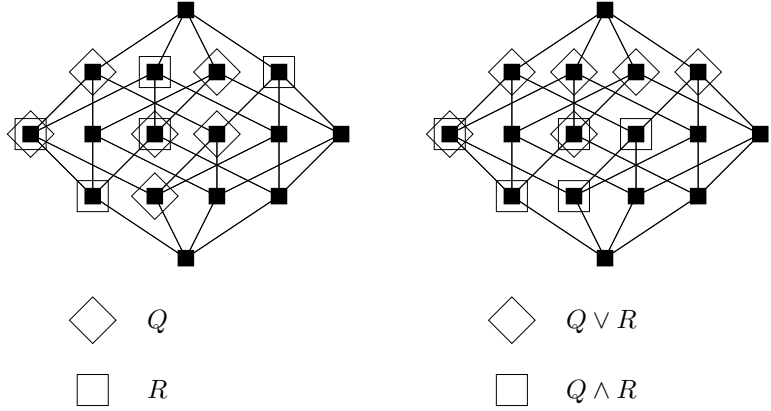}
\end{center}
\caption{An example of the operation $(Q,R) \to (Q\vee R, Q \wedge
R)$ inside the boolean lattice $B_4$.}\label{fig:sch33}
\end{figure}

In Figure \ref{fig:sch33} an example of the operations $\vee$ and
$\wedge$ for two convex subposets of the Boolean lattice $B_4$ is
given. One can easily check that $Q$ and $R$ are indeed convex and
$Q \vee R$ and $Q \wedge R$ are indeed formed according to the rules
above.

\begin{figure}
\begin{center}
\epsfig{file=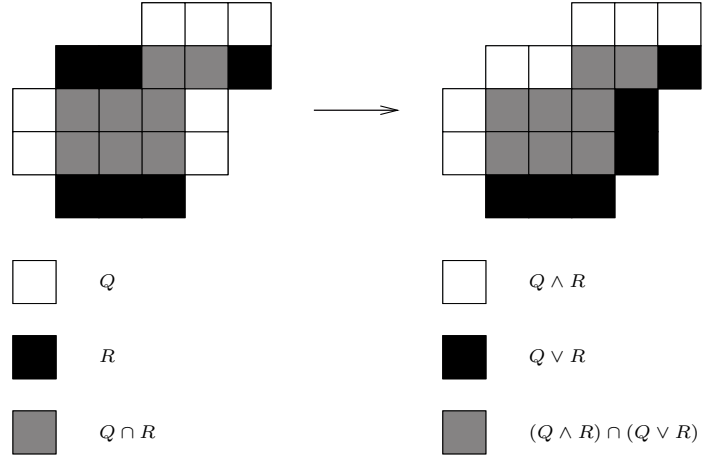}
\end{center}
\caption{An example of the operation $(Q,R) \to (Q\vee R, Q \wedge
R)$ for skew Young diagrams.}\label{fig:sch2}
\end{figure}

Recall that cells in a skew Young diagram form a partially ordered
set, where each cell is covered by the neighboring cell on the right
and the neighboring cell below.  Figure~\ref{fig:sch2} gives an
example of the operations $\vee$ and $\wedge$ for $Q =
(6,5,5,5)/(3,3)$ and $R=(6,6,4,4,4)/(6,1,1,1,1)$, treated as
subposets of the poset $\N^2$ of boxes in the plane (see
Section~\ref{sec:sym}).

Recall that if $A$ and $B$ are sets then $A
\backslash B = \{a \in A \mid a \notin B\}$ denotes the {\it set
difference}.

\begin{lemma}
The subposets $Q \wedge R$ and $Q \vee R$ are both convex
subposets of $P$. We have $(Q \wedge R) \cup (Q \vee R) = Q \cup
R$ and $(Q \wedge R) \cap (Q \vee R) = Q \cap R$.
\end{lemma}
\begin{proof}
We show that $Q \wedge R$ is convex; the statement for $Q \vee R$
follows similarly.  Suppose $s < t$ lie in $Q \wedge R$ and $s < r <
t$ for some $r \in P$ but $r \notin Q \wedge R$. Then $t \in Q
\wedge R$ implies either $t < Q$ or $t \in Q$.  Since $r < t$, we
have either $r < Q$ or $r \in Q$.

If $r<Q$ then $s < Q$ and therefore $s \in R$. Thus either $r \in R$
or $r > R$. If $r \in R$ then since $r < Q$, we get $r \in Q \wedge
R$, obtaining a contradiction. If $r > R$, then $t>r$ implies $t >
R$ which contradicts $t \in Q \wedge R$.

If $r \in Q$ then $r \notin Q \wedge R$ implies $r > R$, and we
proceed as above.

The second statement of the lemma is straightforward.
\end{proof}

Note that the operations $\wedge$ and $\vee$ are stable so that
$(Q \wedge R) \wedge (Q \vee R) = Q \wedge R$ and $(Q \wedge R)
\vee (Q \vee R) = Q \vee R$.

Let $\omega$ be a $(Q,O)$-tableau and $\sigma$ be an
$(R,O)$-tableau. We now describe how to construct a $(Q \wedge R,
O)$-tableau $\omega \wedge \sigma$ and a $(Q \vee R, O)$-tableau
$\omega \vee \sigma$. Define a subset of $Q \cap R$, depending on
$\omega$ and $\sigma$, by
\begin{align*}
(Q \cap R)^+ &= \{x \in Q \cap R \mid \omega(x) < \sigma(x)\}.
\end{align*}
We give $(Q \cap R)^+$ the structure of a graph by inducing from
the Hasse diagram of $Q\cap R$.

Let $\bd(R) = \{ x \in Q \cap R \mid x \gtrdot y \;\mbox{for some
$y \in R \backslash Q$} \} $ be the ``lower boundary'' of $Q \cap
R$ which touches elements in $R$.  Let $\bd(R)^+ \subset (Q \cap
R)^+ $ be the union of the connected components of $(Q \cap R)^+$
which contain an element of $\bd(R)$. Similarly, let $\bd(Q) = \{
x \in Q \cap R \mid x \lessdot y \;\mbox{for some $y \in Q
\backslash R$} \} $ be the ``upper boundary'' of $Q \cap R$ which
touches elements in $Q$.  Let $\bd(Q)^+ \subset (Q \cap R)^+ $ be
the union of the connected components of $(Q \cap R)^+$ which
contain an element of $\bd(Q)$. The elements in $ \bd(Q)^+ \cup
\bd(R)^+$ are amongst the cells that we might ``transfer''.

Let $S \subset Q \cap R$.  Define $(\omega \wedge \sigma)_S: Q
\wedge R \rightarrow \P$ by
\begin{equation}
\label{eq:S1} (\omega \wedge \sigma)_S(x) = \begin{cases}
\sigma(x) & \mbox{if $x \in R \backslash Q$ or $x \in S$,} \\
\omega(x) & \mbox{otherwise.}
\end{cases}
\end{equation}
Similarly, define $(\omega \vee \sigma)_S: Q \vee R \rightarrow \P$ by
\begin{equation}
\label{eq:S2} (\omega \vee \sigma)_S(x) = \begin{cases}
\omega(x) & \mbox{if $x \in Q \backslash R$ or $x \in S$,} \\
\sigma(x) & \mbox{otherwise.}
\end{cases}
\end{equation}

One checks directly that $\wt(\sigma) + \wt(\omega) = \wt((\omega
\wedge \sigma)_S) + \wt((\omega \vee \sigma)_S)$.

\begin{proposition}
\label{prop:respect} Let $(P,O)$ be a $\T$-labelled poset, $Q$ and
$R$ be convex subposets of $P$, and $\omega$ and $\sigma$ be a
$(Q,O)$-tableau and an $(R,O)$-tableau respectively.  Let $S^* :=
\bd(Q)^+ \cup \bd(R)^+$. Then both $(\omega \wedge \sigma)_{S^*}$
and $(\omega \vee \sigma)_{S^*}$ respect $O$.
\end{proposition}

\begin{proof}
We check this for $(\omega \wedge \sigma)_{S^*}$ and the claim for
$(\omega \vee \sigma)_{S^*}$ follows from symmetry. Let $s \lessdot t$
be a covering relation in $Q \wedge R$. Since
$\sigma$ and $\omega$ are assumed to respect $O$, we need only
check the conditions when $(\omega \wedge \sigma)_{S^*}(s) =
\omega(s) (\neq \sigma(s))$ and $(\omega \wedge \sigma)_{S^*}(t) =
\sigma(t) (\neq \omega(t))$; or when $(\omega \wedge
\sigma)_{S^*}(s) = \sigma(s) (\neq \omega(s))$ and $(\omega \wedge
\sigma)_{S^*}(t) = \omega(t) (\neq \sigma(t))$.

In the first case, we must have $s \in Q$ and $t \in R$.  If $t
\in R$ but $t \notin Q$ then by the definition of $Q \wedge R$ we
must have $t < Q $ and so $t < t'$ for some $t' \in Q$. This is
impossible since $Q$ is convex.  Thus $t \in Q \cap R$ and so $t
\in S^*$.  We compute that $\omega(s) \leq O(s,t)(\omega(t)) \leq
O(s,t)(\sigma(t))$ since $\omega(t) < \sigma(t)$ and $O(s,t)$ is
weakly increasing.

In the second case, we must have $s \in R$ and $t \in Q$.  By the
definition of $Q \wedge R$ we must have $t \in R$ as well.  So $t
\in Q \cap R$ but $t \notin S^*$ which means that $\omega(t) >
\sigma(t)$. Thus $\sigma(s) \leq O(s,t)(\sigma(t)) \leq
O(s,t)(\omega(t))$ and $(\omega \wedge \sigma)_{S^*}$ respects $O$ here.
\end{proof}

For each $(\omega,\sigma)$, we say a subset $S \subseteq S^*$ is
{\it transferrable} if both $(\omega \wedge \sigma)_{S}$ and
$(\omega \vee \sigma)_{S}$ respect $O$.

\begin{lemma}
\label{lem:trans} If $S'$ and $S''$ are both transferrable then so
is $S' \cap S''$.
\end{lemma}

\begin{proof}
Let $s \lessdot t$ be a covering relation of $Q \cup R$. Then the
pair $((\omega \wedge \sigma)_{S' \cap S''}(s), (\omega \wedge
\sigma)_{S' \cap S''}(t))$ coincides with one of the four pairs
$(\omega(s),\omega(t))$, $(\sigma(s),\sigma(t))$, $((\omega \wedge
\sigma)_{S'}(s), (\omega \wedge \sigma)_{S'}(t))$ or $((\omega
\wedge \sigma)_{S''}(s), (\omega \wedge \sigma)_{S''}(t))$,
depending on the memberships and non-memberships of $s$, $t$ in $S'$, $S''$.
Since all these pairs are compatible with $O$, so is the pair $((\omega
\wedge \sigma)_{S' \cap S''}(s), (\omega \wedge \sigma)_{S' \cap
S''}(t))$.  The same argument applies for $((\omega \vee \sigma)_{S'
\cap S''}(s), (\omega \vee \sigma)_{S' \cap S''}(t))$.
\end{proof}

Lemma~\ref{lem:trans} implies that there exists a unique smallest
transferrable subset $S^\diamond \subseteq S^*$.  The set
$S^\diamond$ is the key subset used in the proof of the Cell Transfer
Theorem below.  It is the set of ``transferred cells''.

Define $\eta: \mathcal{A}(Q ,O)
\times \mathcal{A}(R ,O) \rightarrow \mathcal{A}(Q \wedge R ,O)
\times \mathcal{A}(Q \vee R ,O)$ by
\[
(\omega,\sigma) \longmapsto ((\omega \wedge
\sigma)_{S^\diamond},(\omega \vee \sigma)_{S^\diamond}).
\]
Note that $S^\diamond$ depends on $\omega$ and $\sigma$, though we
have suppressed the dependence from the notation.

We call the map $\eta$ the {\it {cell transfer procedure}}.  This
name comes from our motivating example, where the elements of the
poset are the cells of a Young diagram $\lambda$.  For convenience,
in the following proof, we call elements of any poset $P$ {\it
{cells}}. We say that a cell $s \in P$ is {\emph{transferred}} if $s
\in S^\diamond$.


\begin{lemma}
\label{lem:inj} The map $\eta$ is injective.
\end{lemma}

\begin{proof}
Given $(\alpha,\beta)
\in \eta(\mathcal{A}(Q ,O) \times \mathcal{A}(R ,O))$, we show how
to recover $\omega$ and $\sigma$.  As before, for a subset $S
\subset Q \cap R$, define $\omega_S= \omega(\alpha,\beta)_S: Q
\rightarrow \P$ by
\[
\omega_S(x)  = \begin{cases} \beta(x) & \mbox{if $x \in (Q
\backslash R) \cap (Q \vee R)$ or $x \in S$,} \\
\alpha(x) & \mbox{otherwise.}
\end{cases}
\]
And define $\sigma_S = \sigma(\alpha,\beta)_S: R \rightarrow \P$
by
\[
\sigma_S(x) = \begin{cases} \alpha(x) & \mbox{if $x \in (R
\backslash Q) \cap (Q \wedge R)$ or $x \in S$,} \\
\beta(x) & \mbox{otherwise.}
\end{cases}
\]
Note that if $(\alpha,\beta) = ((\omega \wedge
\sigma)_{S^\diamond},(\omega \vee \sigma)_{S^\diamond}))$ then
$\omega = \omega_{S^\diamond}$ and $\sigma = \sigma_{S^\diamond}$.
Let $S^\square \subset Q \cap R$ be the unique smallest subset
such that $\omega_{S^\square}$ and $\sigma_{S^\square}$ both
respect $O$.  Since we have assumed that $(\alpha,\beta) \in
\eta(\mathcal{A}(Q ,O) \times \mathcal{A}(R ,O))$, such a
$S^\square$ must exist.  (As before the intersection of two
transferrable subsets with respect to $(\alpha,\beta)$ is
transferrable.)

We now show that if $(\alpha,\beta) = ((\omega \wedge
\sigma)_{S^\diamond},(\omega \vee \sigma)_{S^\diamond}))$ then
$S^\square = S^\diamond$.  We know that $S^\square \subset
S^\diamond$ from the previous paragraph.  Let $C \subset S^\diamond
\backslash S^\square$ be a connected component of $S^\diamond
\backslash S^\square$, viewed as an induced subgraph of the Hasse
diagram of $P$.  We claim that $S^\diamond \backslash C$ is a
transferrable set for $(\omega,\sigma)$; this means that changing
$\alpha|_C$ to $\omega|_C$ and $\beta|_C$ to $\sigma|_C$ gives a
pair in $\mathcal{A}(Q \wedge R ,O) \times \mathcal{A}(Q \vee R
,O)$. Suppose first that $c \in C$ and $s \in S^\square$ is so that
$c \lessdot s$. By the definition of $S^\square$, we must have
$\alpha(c) \leq O(c,s)(\beta(s))$ and $\beta(c) \leq
O(c,s)(\alpha(s))$.  Now suppose that $c \in C$ and $s \in Q
\backslash R$ such that $c \lessdot s$.  Then we must have
$O(c,s)(\omega(s)) = O(c,s)(\beta(s)) \geq \alpha(c) = \sigma(c)$.
Similar conclusions hold for $c \gtrdot s$.  Thus we have checked
that $S^\diamond \backslash C$ is a transferrable set for
$(\omega,\sigma)$, which is impossible by definition of
$S^\diamond$: it is the minimal set with this property.  Therefore
the set $S^\diamond \backslash S^\square$ is empty and thus
$S^\diamond = S^\square$.

Thus the map $\mu: \eta(\mathcal{A}(Q ,O) \times \mathcal{A}(R ,O))
\to \mathcal{A}(Q ,O) \times \mathcal{A}(R ,O)$ given by
$$
\mu: (\alpha,\beta) \longmapsto (\omega(\alpha,\beta)_{S^\square},
\sigma(\alpha,\beta)_{S^\square})$$ is inverse to $\eta$. This shows
that the map $(\omega,\sigma) \mapsto ((\omega \wedge
\sigma)_{S^\diamond},(\omega \vee \sigma)_{S^\diamond}))$ is
injective, completing the proof.
\end{proof}

We say that a map between pairs of tableaux
{\emph{weight-preserving}} if the multiset of their values over all
$s \in P$ is not changed by the map.

\begin{theorem}[Cell Transfer Theorem]
\label{thm:celltransfer} The difference \[K_{Q \wedge R,O}K_{Q\vee
R,O} - K_{Q,O}K_{R,O}\] is monomial-positive.
\end{theorem}

\begin{proof}
The map
$$\eta: \mathcal{A}(Q ,O) \times \mathcal{A}(R ,O) \longrightarrow
\mathcal{A}(Q \wedge R ,O) \times \mathcal{A}(Q \vee R ,O)$$ defined
above is weight-preserving.  Indeed, for each element $s \in
P \cup Q$ we have $\{\omega(s), \sigma(s)\} = \{(\omega \wedge
\sigma)_{S^\diamond}(s), (\omega \vee \sigma)_{S^\diamond}(s)\}$ as
multisets, where the value of a tableau is zero outside of its range
of definition. Then since the map $\eta$ is injective and since
$K_{Q \wedge R,O} K_{Q \vee R,O}$ and $K_{Q,O}K_{R,O}$ are the
generating functions of common weights of pairs of the tableaux of
corresponding shapes, the statement follows.
\end{proof}


\begin{figure}
\begin{center}
\epsfig{file=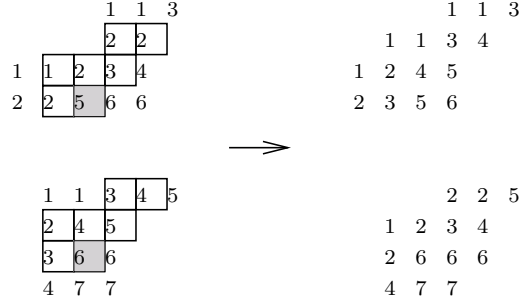}
\end{center}
\caption{An example of the map $\eta$ applied to a pair of semistandard Young
tableaux.  The cells in $S^*$ are marked, and the unique cell in $S^{*}/S^{\diamond}$ is shaded.}\label{fig:ss}
\end{figure}

In Figure \ref{fig:ss} an example of the cell transfer injection $\eta$ is given for a pair of
tableaux with the shapes $Q = (6,5,5,5)/(3,3)$ and $R=(6,6,4,4,4)/(6,1,1,1,1)$ that were shown in Figure \ref{fig:sch2}. Note that
there is one cell contained in $S^*$ but not $S^{\diamond}$: the cell
labeled $5$ in $Q$ and $6$ in $R$.

Note that $(\omega,\sigma) \mapsto ((\omega \wedge
\sigma)_{S^*},(\omega \vee \sigma)_{S^*})$ also defines a
weight-preserving map $\eta^*: \mathcal{A}(Q ,O) \times
\mathcal{A}(R ,O) \rightarrow \mathcal{A}(Q \wedge R ,O) \times
\mathcal{A}(Q \vee R ,O)$.  Unfortunately, $\eta^*$ is not always
injective.

Suppose $P$ is a locally-finite poset with a unique minimal element.
Let $J(P)$ be the lattice of finite order ideals of $P$;
see~\cite{EC2}. If $I, J \in J(P)$ then the subposets $I \wedge J$
and $I \vee J$ of $P$ defined in (\ref{eq:wedge}) and (\ref{eq:vee})
are finite order ideals of $P$ and agree with the the meet
$\wedge_{J(P)}$ and join $\vee_{J(P)}$ of $I$ and $J$ respectively
within $J(P)$.  In this case the operations $\vee$ and $\wedge$ are
commutative.

Now let $P$ be any poset.  By defining $Q \wedge' R = \{s \in R \mid
s < Q \} \cup \{s \in Q \mid s \in R \;\text{or}\; s < R\}$ and $Q
\vee' R = \{s \in Q \mid s \sim R \;\text{or}\; s >R\} \cup \{s \in
R \mid s \sim Q \;\text{or}\; s > Q\}$, the order ideals $I \wedge'
J = I \wedge_{J(P)} J$ and $I \vee' J = I \vee_{J(P)} J$ agree with
the meet and join in $J(P)$ even when $P$ does not contain a minimal
element.


\begin{proposition}
Let $P$ be a locally-finite poset, let $I, J$ be elements of  $
J(P)$, and let $O$ be a $\mathbb T$-labeling of $P$.  Then the
generating function
\[
K_{I \wedge_{J(P)} J, O} K_{I \vee_{J(P)} J, O} - K_{I,O}K_{J,O}
\]
is monomial-positive.
\end{proposition}
\begin{proof}
Let $(\omega,\sigma) \in \mathcal{A}(Q ,O) \times \mathcal{A}(R
,O)$.  Replacing $Q \wedge R$ by $Q\wedge' R$ and $Q \vee R$ by
$Q\vee'R$ in (\ref{eq:S1}) and (\ref{eq:S2}) we can define $(\omega
\vee' \sigma)_S$ and $(\omega \wedge' \sigma)_S$.

The conclusion of Proposition~\ref{prop:respect} holds with $\wedge$ and
$\vee$ replaced by $\wedge'$ and $\vee'$.  This is because the set $C$ of elements of $Q \wedge R$
not belonging to $Q\wedge' R$ are exactly the elements $s \in Q$ which are incomparable with elements
of $R$.  These elements belong instead to $Q \vee' R$.  Since the cells in $C$ are incomparable with elements of $R$, they are in particular never compared with the elements of $Q \cap R$ in the proof
of Proposition~\ref{prop:respect}.  Thus to show that $(\omega
\vee' \sigma)_{S^*}$ and $(\omega \wedge' \sigma)_{S^*}$ respect $O$ the same set of inequalities needs to be verified as in the original Proposition~\ref{prop:respect}.

Lemma~\ref{lem:trans} also holds with a verbatim proof if we define
$S \subset S^*$ to be transferrable if both $(\omega \vee'
\sigma)_{S}$ and $(\omega \wedge' \sigma)_{S}$ respect $O$.

Finally, using the same definition (following Lemma~\ref{lem:trans})
of the set $S^\diamond$, we can obtain a map $\eta': \mathcal{A}(Q
,O) \times \mathcal{A}(R ,O) \rightarrow \mathcal{A}(Q \wedge' R ,O)
\times \mathcal{A}(Q\vee'R ,O)$ analogous to $\eta$. By the modified
versions of Proposition~\ref{prop:respect} and Lemma~\ref{lem:trans}
the image of $\eta'$ consists of pairs of $O$-compatible labelings.
The map $\eta'$ is also injective: the proof of Lemma~\ref{lem:inj}
remains valid since the cells in $C$ are incomparable with the
elements in $S^\diamond$, and the calculations in the proof of
Lemma~\ref{lem:inj} always involve some element of $S^\diamond$.

Now the proof of Theorem~\ref{thm:celltransfer} can be modified by replacing $\wedge$, $\vee$ and $\eta$
with $\wedge'$, $\vee'$ and $\eta'$ to obtain the claimed statement.

\end{proof}

\section{Symmetric and Quasisymmetric functions}
\label{sec:sym} We refer to~\cite{EC2} for more details of the
material in this section.

Let $n$ be a positive integer.  A {\it composition} of $n$ is a
sequence $\alpha = (\alpha_1,\alpha_2,\ldots,\alpha_k)$ of
positive integers such that $\alpha_1 + \alpha_2 + \cdots +
\alpha_k = n$.  We write $|\alpha| = n$.  If in addition $\alpha_1
\geq \alpha_2 \geq \cdots \geq \alpha_k$ then we say that $\alpha$
is a {\it partition} of $n$.  If $\lambda$ is a partition then
$\lambda'$ denotes the {\it conjugate} partition.  Let $l(\ll)$
denote the number of (non-zero) parts of $\lambda$.


A formal power series $f = f(x) \in \Q[[x_1,x_2,\ldots]]$ with
bounded degree is called {\it quasisymmetric} if for any
$a_1,a_2,\ldots,a_k \in \P$ we have
\[
\left[x_{i_1}^{a_1}\cdots x_{i_k}^{a_k} \right]f =
\left[x_{j_1}^{a_1}\cdots x_{j_k}^{a_k} \right]f
\]
whenever $i_1 < \cdots < i_k$ and $j_1 <\cdots < j_k$.  Here
$[x^\alpha]f$ denotes the coefficient of $x^\alpha$ in $f$. Denote
by $\Qsym \subset  \Q[[x_1,x_2,\ldots]]$ the space (in fact
algebra) of quasisymmetric functions.

Let $\alpha$ be a composition. Then the {\it monomial
quasisymmetric function} $M_\alpha$ is given by
\[
M_\alpha = \sum_{i_1 < \cdots <i_k} x_{i_1}^{\alpha_k} \cdots
x_{i_k}^{\alpha_k}.
\]
The set of monomial quasisymmetric functions form a basis of
$\Qsym$. Another basis is given by the {\it fundamental quasi-symmetric functions} $L_\alpha$ defined as follows:
\[
L_\alpha = \sum_{\beta \leq \alpha} M_\beta,
\]
where for two compositions $\alpha,\beta$ we have $\beta \leq
\alpha$ if and only if $\beta$ is a refinement of $\alpha$.

Define (as in Example~\ref{ex:schur}) two functions
$f^{\weak},f^{\strict}:\P \rightarrow \N \cup \{\infty\}$ by
$f^{\weak}(n) = n$ and $f^{\strict}(n) = n-1$.

\begin{proposition}
Let $(P,O)$ be a finite $\T$-labelled poset.  Suppose $$O(s,t) \in
\{f^{\weak},f^{\strict}\}$$ for each covering relation $s \lessdot
t$.  Then $K_{P,O}(x)$ is a quasi-symmetric function.
\end{proposition}

\begin{example}
Let $P_n$ denote the chain with $n+1$ elements as in Example~\ref{ex:chain}.
\end{example}

A $\T$-labelled poset satisfying the conditions of the proposition
is called {\it oriented} in~\cite{McN}.  Stanley's
$(P,\omega)$-partitions are special cases of $(P,O)$-tableaux, for
such posets.  If $f \in \Qsym$ then $f$ is $m$-positive if and only
if is a non-negative linear combination of the $M_\alpha$.

\medskip
A formal power series $f = f(x) \in \Q[[x_1,x_2,\ldots]]$ with
bounded degree is called {\it symmetric} if for any
$a_1,a_2,\ldots,a_k \in \P$ we have
\[
\left[x_{i_1}^{a_1}\cdots x_{i_k}^{a_k} \right]f =
\left[x_{j_1}^{a_1}\cdots x_{j_k}^{a_k} \right]f
\]
whenever $i_1, \ldots, i_k$ are all distinct and $j_1, \ldots,j_k$
are all distinct.  Denote by $\Lambda \subset
\Q[[x_1,x_2,\ldots]]$ the algebra of symmetric functions.  Every
symmetric function is quasisymmetric.

Given $\ll = (\ll_1,\ll_2,\ldots)$, the monomial symmetric
functions $m_\ll$ is given by
\[
m_\ll(x) = \sum_\alpha x_1^{\alpha_1}\cdots x_k^{\alpha_k}
\]
where the sum is over all distinct permutations $\alpha$ of the
entries of the (infinite) vector $(\ll_1,\ll_2,\ldots)$.  As $\ll$
ranges over all partitions, the $m_\ll$ form a basis of $\Lambda$.
If $f \in \Lambda$ then $f$ is {\it monomial-positive} if and only
if is a non-negative linear combination of the monomial symmetric
functions.

Let $\lambda$ be a partition.  Recall that a {\it semistandard Young
tableau} with shape $\lambda$ is a filling of the squares of the
Young diagram of $\lambda$ with positive integers so that the rows
are weakly increasing and the columns are strictly increasing.  The
{\it Schur function} $s_\lambda(x_1,x_2,\ldots)$ is the following
generating function:
$$
s_\lambda(x_1,x_2,\ldots) = \sum_T x_1^{\#\text{1's in
$T$}}x_2^{\#\text{2's in $T$}}\cdots \, ,
$$
where the summation is over all semistandard Young tableaux of shape
$\lambda$.  More generally one defines the {\it skew Schur
functions} $s_{\lambda/\mu}(x_1,x_2,\ldots)$ in the same manner.  The
Schur functions $s_\lambda$ form a basis of $\Lambda$ as $\lambda$
varies over all partitions. If $f \in \Lambda$ is a non-negative
linear combination of Schur functions then we call $f$ {\it
Schur-positive}.

In the following, we will consider all (skew) shapes as convex subposets of
the poset $\N^2$ of boxes in the plane with partial order $(i,j) \geq (i',j')$ if and only if $i \geq i'$ and $j \geq j'$.
Thus the $i$-th row of a shape $\ll$ has cells with coordinates
$\{(i,1),(i,2),\ldots,(i,\ll_i)\}$.  Note that we use the word {\it shape} to denote a specific such subposet, which may
still have multiple representations of the form $\ll/\mu$.  For example,
 $(6,5,5,5)/(3,3)$ and $(6,5,5,5,1)/(3,3,1)$ represent the same shape but $(6,5,5,5)/(3,3) \neq (7,6,6,6)/(4,4,1,1)$.

\begin{theorem} \label{thm:schur}
The symmetric function $s_{\mu \wedge \lambda}s_{\mu \vee
\lambda}-s_{\mu}s_{\lambda}$  is monomial-positive.
\end{theorem}

\begin{proof}
This follows immediately from Theorem \ref{thm:celltransfer}.
\end{proof}
\noindent

Let $\ll =
(\ll_1,\ldots,\ll_k)$, $ \mu = (\mu_1,\ldots,\mu_k)$, $\nu =
(\nu_1,\ldots,\nu_k)$ and $\rho = (\rho_1,\ldots,\rho_k)$ be four partitions such that
$\mu \subset \ll$ and $\rho \subset \nu$.  Define the shapes
\[
\max(\ll/\mu,\nu/\rho) :=
(\max(\ll_1,\nu_1),\ldots,\max(\ll_k,\nu_k))/(\max(\mu_1,\rho_1),\ldots,\max(\mu_k,\rho_k))
\]
and
\[
\min(\ll/\mu,\nu/\rho) :=
(\min(\ll_1,\nu_1),\ldots,\min(\ll_k,\nu_k))/(\min(\mu_1,\rho_1),\ldots,\min(\mu_k,\rho_k)).
\]
Note that $\max$ and $\min$ depends on all four partitions $\ll,\mu,\nu,\rho$ and not just the
shapes $\ll/\mu$ and $\nu/\rho$.
These shapes $\max(\ll/\mu,\nu/\rho)$ and $\min(\ll/\mu,\nu/\rho)$ are nearly but not always the same as $\ll/\mu \vee
\nu/\rho$ and $\ll/\mu \wedge \nu/\rho$ respectively.  This is
because we may have $\ll_i = \mu_i = a$ for some $i$ and then the
shape $\ll/\mu$ does not depend on the exact value of $a$.  However, as one can see from the definitions above, the shapes $\max(\ll/\mu,\nu/\rho) $ and $\min(\ll/\mu,\nu/\rho)$ do depend on the choice of $a$.

Fix four partitions $\ll,\mu,\nu,\rho$ such that $\mu \subset \ll$
and $\rho \subset \rho$.  Let $V(\ll/\mu,\nu/\rho) \subset (\ll/\mu
\cup \nu/\rho)$ denote the set of cells for which
$\min(\ll/\mu,\nu/\rho)$ and $\ll/\mu \wedge \nu/\rho$ differ.  Here
$\ll/\mu \cup \nu/\rho$ denotes the set theoretic union of the cells
lying in $\ll/\mu$ and $\nu/\rho$. Clearly each cell $s \in
V(\ll/\mu,\nu/\rho)$ lies in only one of $\ll/\mu$ or $\nu/\rho$.

\begin{lemma}
\label{lem:minmaxdifference} Let $s \in V(\ll/\mu,\nu/\rho)$.  If $s
\in \ll/\mu$ then $s$ is incomparable to all squares in $\nu/\rho$.
If $s \in \nu/\rho$ then $s$ is incomparable to all squares in
$\ll/\mu$.
\end{lemma}

\begin{proof}
Let $V_i \subset V(\ll/\mu,\nu/\rho)$ denote the set of cells in the
$i$-th row for which $(\ll/\mu \vee \nu/\rho,\ll/\mu \wedge
\nu/\rho)$ and $(\max(\ll/\mu,\nu/\rho),\min(\ll/\mu,\nu/\rho))$
differ.  If $V_i$ is non-empty then either $\ll_i = \mu_i$ or
$\nu_i=\rho_i$ (but not both).

Without loss of generality we assume that $\ll_i=\mu_i = a$
so that $V_i \subset \nu/\rho$.  Let
the leftmost cell in the lowest non-empty row of $\ll/\mu$ above row
$i$ have coordinates $(p,a')$ and let the rightmost cell in the
highest non-empty row of $\ll/\mu$ below row $i$ have coordinates
$(q,a'')$. Then in particular $a'' \leq a \leq a'$.  It is easy to
check, case by case, that $S_i \subset
\{(i,a''+1),(i,a''+1),\ldots,(i,a')\}$. These cells are incomparable
with any cells in $\ll/\mu$.

\end{proof}

\begin{theorem} \label{thm:skewschur}
The symmetric function
$s_{\max(\ll/\mu,\nu/\rho)}s_{\min(\ll/\mu,\nu/\rho)}-s_{\ll/\mu}s_{\nu/\rho}$
is monomial-positive.
\end{theorem}

\begin{proof}
We give an injection from the set of pairs $(U,T)$ of semistandard
tableaux of shape $(\ll/\mu,\nu/\rho)$ to the set of pairs $(U',T')$
of semistandard tableaux with shape
$(\max(\ll/\mu,\nu/\rho),\min(\ll/\mu,\nu/\rho))$. First we apply
the map $\eta$ of Theorem~\ref{thm:celltransfer} to $(U,T)$ to
obtain a pair $(U'',T'')$ of semistandard tableaux of shape
$(\ll/\mu \vee \nu/\rho,\ll/\mu \wedge \nu/\rho)$.

Now define $U'$ by letting it be the unique tableau of shape
$\max(\ll/\mu,\nu/\rho)$ with the same numbers as $U''$ in the boxes
of $\max(\ll/\mu,\nu/\rho) \cap (\ll/\mu \vee \nu/\rho)$ and with
the same numbers as $T''$ in the boxes of $\max(\ll/\mu,\nu/\rho)
\cap V(\ll/\mu,\nu/\rho)$.  Similarly define $T'$.  This process is
reversible so the map $(U,T) \to (U',T')$ is injective.  We claim
that $(U',T')$ is still semistandard, from which the theorem
follows.


The claim follows from Lemma~\ref{lem:minmaxdifference}: if a cell $s \in
V(\ll/\mu,\nu/\rho)$ was originally contained in $\nu/\rho$
(respectively $\ll/\mu$) then there is no cell of $\ll/\mu$
(respectively $\nu/\rho$) adjacent to it.  Thus checking the
inequalities assuring semistandard-ness for a cell $s \in V$ is
trivial: the inequality was satisfied in the semistandard tableau
$U$ or $T$.

\end{proof}

\begin{conjecture}
\label{conj:LPP} The symmetric function
$s_{\max(\ll/\mu,\nu/\rho)}s_{\min(\ll/\mu,\nu/\rho)}-s_{\ll/\mu}s_{\nu/\rho}$
is Schur-positive.
\end{conjecture}

This conjecture is proved in joint work~\cite{LPP} with Alex
Postnikov.  The proof relies ultimately on some deep results in
representation theory.

A combinatorial proof of the weaker statement that this same
expression is positive in terms of fundamental quasisymmetric
functions is given in \cite{LP}.

\section{Cell transfer as an algorithm}

Let $(P, O)$ be a $\mathbb T$-labelled poset. We now give an algorithmic description of cell transfer.  Let $Q$
and $R$ be two finite convex subposets of $P$. We construct
step-by-step an injection
$$\eta: \mathcal{A}(Q ,O) \times \mathcal{A}(R ,O) \longrightarrow
\mathcal{A}(Q \wedge R ,O) \times \mathcal{A}(Q \vee R ,O)$$ which
is weight-preserving. Let $\omega$ be a $(Q,O)$-tableau and
$\sigma$ be an $(R,O)$-tableau.
Let us recursively define $\bar \omega: Q \wedge R \rightarrow
\mathbb P$ and $\bar \sigma: Q \vee R \rightarrow \mathbb P$ as
follows.

\begin{enumerate}
\item Define $\bar \omega: Q \wedge R \rightarrow \mathbb P$ and $\bar \sigma: Q \vee R \rightarrow \mathbb P$ as follows:
\begin{align*}
\bar \omega(s) &=
\begin{cases}
\omega(s) & \text{if $s \in Q$,}\\
\sigma(s) & \text{if $s \in R/Q$.}\\
\end{cases}\\
\bar \sigma(s) &=
\begin{cases}
\sigma(s) & \text{if $s \in R$,}\\
\omega(s) & \text{if $s \in Q/R$.}\\
\end{cases}
\end{align*} Note that $\bar \omega$ and $\bar \sigma$ do not necessarily
respect $O$.  Indeed, the parts of $\omega$ and $\sigma$ which we
glued together might not agree with each other, i.e., a covering
relation $s \lessdot t$ might fail to respect $O$, where the label
of one of $s, t$ comes from $\omega$, and that of the other from
$\sigma$.
\item
We say that we {\it transfer} a cell $s \in Q \cap R$ when we swap
the values at $s$ of $\bar \omega$ and $\bar \sigma$. We say that
a cell $s$ in $Q \cap R$ is {\it {critical}} if one of the
following condition holds
\begin{enumerate}
\item for some $t \in R$ and  $s \gtrdot t$ we have $O(s, t)(\bar \omega(s)) < \bar \omega(t)$,
\item for some $t \in Q \cap R$ and $t \gtrdot s$ we have $O(s, t)(\bar \sigma(t)) < \bar \sigma(s)$,
\item for some $t \in Q$ and $t \gtrdot s$ we have $O(s, t)(\bar \sigma(t)) < \bar \sigma(s)$,
\item for some $t \in Q \cap R$ and $s \gtrdot t$ we have $O(s, t)(\bar \omega(s)) < \bar \omega(t)$,
\end{enumerate}
and $s$ was not transferred in a previous iteration.  We now
transfer all critical cells if there are any.
\item Repeat step (2) until no critical cells are transferred.
\end{enumerate}

\begin{theorem}[Cell Transfer Algorithm]
The algorithm described above terminates in a finite number of
steps. The resulting maps $\bar \omega$ and $\bar \sigma$ are
$(P,O)$-tableaux and coincide with $(\omega \wedge
\sigma)_{S^\diamond}$ and $(\omega \vee \sigma)_{S^\diamond}$
defined in the proof of Theorem~\ref{thm:celltransfer}.
\end{theorem}

\begin{proof}
As for the first claim, there is a finite number of cells
in $Q \cap R$ and each gets transferred at most once, thus the
process terminates.

We say that an edge $a \lessdot b$ in the Hasse diagram of $P$
respects $O$ if $\bar \omega(a) \leq O(a,b) (\bar \omega(b))$ and
$\bar \sigma(a) \leq O(a,b) (\bar \sigma(b))$, whenever these
inequalities make sense.  Note that a cell $s$ is critical only if
the cell $t$ (from the definition of a critical cell) was
transferred in previous iteration of step (2), or if it is the
first iteration of step (2) and $t$ belongs to $\{s \in R \mid s <
Q\}$ or $\{s \in Q \mid s > R\}$ -- the parts which were ``glued''
in step (1). Indeed, if $s,t$ have both not been transferred then
$s \lessdot t$ must respect $O$ since $\omega$ and $\sigma$ were
$(P,O)$-tableaux to begin with. Similarly, two cells $s \lessdot
t$ which have both been transferred must also respect $O$.

We thus see that after the second step every edge between $\{s \in
R \mid s < Q\}$ and $Q \cap R$, as well as between $\{s \in Q \mid
s > R\}$ and $Q \cap R$ respects $O$.  After the algorithm
terminates every edge in $Q \cap R$ must respect $O$, since if
there exists an edge $s \lessdot t$ which does not then one of $s$
and $t$ must have already been transferred, and the other has not
been transferred and thus is critical.  This contradicts the
termination condition of the algorithm.  Therefore, the only
possible edges which might fail to respect $O$ are the ones
between $\{s \in Q  \mid s < R\}$ and $Q \cap R$, and the ones
between $\{s \in R \mid s > Q\}$ and $Q \cap R$. However, it is
easy to see that during the whole process values of $\bar \omega$
on $Q \cap R$ are increasing, and therefore cannot be not large
enough for values of $\bar \omega$ on $\{s \in Q \mid s < R\}$.
Similarly, values of $\bar \sigma$ on $Q \cap R$ are decreasing
and cannot be too large for values of $\bar \sigma$ on $\{s \in R
\mid s > Q\}$. Thus, we do obtain two $(P,O)$-tableaux $\bar
\omega$ and $\bar \sigma$.

Let $\bar S \subset Q \cap R$ be the set of cells we transferred
during the algorithm. The fact that values of $\bar \omega$ on $Q
\cap R$ increase and the values of $\bar \sigma$ on $Q \cap R$
decrease implies that $\bar S$ is contained in $S^*$ (as defined
in the proof of Theorem~\ref{thm:celltransfer}). We claim that all
transferrable sets contain $\bar S$. Indeed, in each iteration we
transfer only those cells that must be transferred in order for
the result to respect $O$. On the other hand, as shown above the
set $\bar S$ is transferrable itself.  Thus, it is exactly the set
$S^\diamond$ -- the minimal transferrable set. This completes the
proof of the theorem.
\end{proof}

The algorithmic description above provides another way to verify
injectivity of $\eta$.  Let $\bar \omega: Q \wedge R \rightarrow
\mathbb P$ and $\bar \sigma: Q \vee R \rightarrow \mathbb P$ be in
the image of $\eta$.  Then one can define maps $\omega': Q
\rightarrow \mathbb P$ and $\sigma': R \rightarrow \P$ by
\begin{align*}
\omega'(s) &=
\begin{cases}
\bar \omega(s) & \text{if $s \in Q \cap(Q \wedge R)$,}\\
\bar \sigma(s) & \text{otherwise.}\\
\end{cases} \\
\sigma(s) &=
\begin{cases}
\bar \sigma(s) & \text{if $s \in R \cap (Q \vee R)$,}\\
\bar \omega(s) & \text{otherwise.}\\
\end{cases}
\end{align*}
We now iterate step (2) of the cell transfer algorithm with
$\omega'$ and $\sigma'$ replacing $\bar \omega$ and $\bar \sigma$.

One can verify that for each step the set of transferred cells is
identical to the corresponding step of the original algorithm for
$\bar \omega$ and $\bar \sigma$.  This produces the inverse of
$\eta$.

\begin{figure}
\begin{center}
\epsfig{file=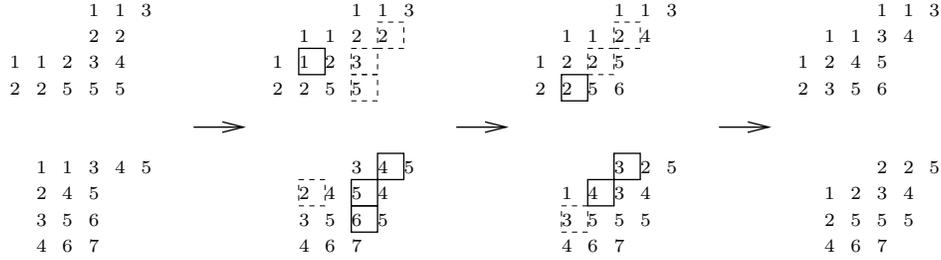}
\end{center}
\caption{An example of the cell transfer
algorithm applied to a pair of semistandard tableaux.  The critical cells are marked with squares.}\label{fig:sch3}
\end{figure}

In Figure~\ref{fig:sch3}, we show the step-by-step application of the cell transfer
algorithm to a pair of semistandard tableaux $(\omega,\sigma)$.


\section{Final Remarks}
The most interesting feature of the cell transfer theorem is that
for the case of Schur functions, the theorem holds with
monomial-positivity replaced by Schur positivity~\cite{LPP}.  We
know of no simple combinatorial explanation of this phenomenon. A
natural question to ask is whether a result similar to
Conjecture~\ref{conj:LPP} holds for other $\T$-labelled posets
$(P,O)$. When $P$ has a minimal element, it seems reasonable to
replace Schur-positivity in the conjecture by positivity in the
generating functions $\{K_{I,O} \mid \mbox{$I$ is an order ideal
of $P$}\}$. However, it is not even clear under what conditions
the functions $\{K_{I,O}\}$ might span the space of functions
$\{K_{Q,O}\mid \mbox{$Q$ is a convex subposet of $P$}\}$ or span
the space of differences of products of such functions.
$\T$-labelled posets satisfying these weaker requirements would
also be worth studying.

Cylindric Schur functions provide another interesting special
case.  It is conjectured in~\cite{Lam} that the cylindric Schur
functions $K_{(P,O)}(x_1,x_2,\cdots)$, where $P$ is a cylindric
skew shape (see Example~\ref{ex:cylindric}) is a positive linear
combination of symmetric functions known as {\it affine Schur
functions} or {\it dual $k$-Schur functions} (see
also~\cite{McN}).  One might speculate that the differences of
products of cylindric Schur functions from
Theorem~\ref{thm:celltransfer} are also affine Schur-positive.

\end{document}